\documentclass[12pt]{article}
\usepackage{amsfonts, amsmath, amssymb,amsthm,bbm,comment,fancyhdr,float,makeidx,mathrsfs,verbatim}

\normalsize

\newcommand{\si}{\sigma}

\newcommand{\seq}{\subseteq}

\newcommand{\ol}{\overline}

\newcommand{\vs}{\vspace*}

\newcommand{\nin}{\noindent}

\newtheorem{mthm}{Theorem}[section]
\newtheorem{mylem}[mthm]{Lemma}
\newtheorem{myprn}[mthm]{Proposition}
\newtheorem{mycor}[mthm]{Corollary}
\newtheorem{mydef}[mthm]{Definition}
\newtheorem{myrem}[mthm]{Remark}
\newtheorem{mycon}[mthm]{Construction}
\newtheorem{myeg} [mthm]{Example}
\newtheorem{myque} [mthm]{Question}
\newtheorem{myalg} [mthm]{Algorithm}
\newenvironment{thm}{\begin{mthm}}{\end{mthm}}
\newenvironment{lem}{\begin{mylem}}{\end{mylem}}
\newenvironment{cor}{\begin{mycor}}{\end{mycor}}
\newenvironment{prop}{\begin{myprn}}{\end{myprn}}
\newenvironment{mdef}{\begin{mydef}\rm}{\end{mydef}}
\newenvironment{rem}{\begin{myrem}\rm}{\end{myrem}}

\newenvironment{ex}{\begin{myeg}\rm}{\end{myeg}}
\newenvironment{prof}{\noindent $Proof.$ \rm}{\hfill $\Box$}

\def \nin {\noindent}

\def \Lemma #1 {\vs{3mm}\nin {\bf Lemma #1} \it}
\def \Prop #1 {\vs{3mm}\nin {\bf Proposition #1} \it}
\def \Th #1 {\vs{3mm}\nin {\bf Theorem #1} \it}
\def \Cor #1 {\vs{3mm}\nin {\bf Corollary #1} \it}
\def \Ex #1 {\vs{3mm}\nin {\bf Example #1} \it}

\def \part #1 {\hfil\break\hglue 12pt {\rm (#1)~}}

\def\fs{\footnotesize}

\title{
\bf\LARGE  Perfect Sets and $f$-Ideals \thanks{This research is supported by the National Natural
Science Foundation of China (Grant No. 11271250). } }
\author{Jin Guo\thanks{ guojinecho@163.com}, \,\, {Tongsuo Wu\thanks{Corresponding author. tswu@sjtu.edu.cn}}\\
 {\small Department of Mathematics, Shanghai Jiaotong University}\\
 Qiong Liu\thanks{ sky200547@126.com}\\
 {\small Department of Mathematics, Shanghai University of Electric Power}
}

\date{}
\begin{document}
\baselineskip=16pt \maketitle

\begin{center}
\begin{minipage}{12cm}

 \vs{3mm}\nin{\small\bf Abstract.} {\fs A square-free monomial ideal $I$ is called an {\it $f$-ideal}, if both $\delta_{\mathcal{F}}(I)$ and $\delta_{\mathcal{N}}(I)$ have the same $f$-vector, where $\delta_{\mathcal{F}}(I)$ ($\delta_{\mathcal{N}}(I)$, respectively) is the facet (Stanley-Reisner, respectively) complex related to $I$. In this paper, we introduce and study perfect subsets of $2^{[n]}$ and use them to characterize  the $f$-ideals of degree $d$. We give a decomposition of $V(n, 2)$ by taking advantage of a correspondence between graphs and sets of square-free monomials of degree $2$, and then give a formula for counting the number of $f$-ideals of degree $2$, where $V(n, 2)$ is the set of $f$-ideals of degree 2 in $K[x_1,\ldots,x_n]$. We also consider the relation between an $f$-ideal and an unmixed monomial ideal.
}


\vs{3mm}\nin {\small Key Words and phrases:} {\small  perfect set; perfect number; $f$-ideal; Stanley-Reisner complex; facet complex }

\vs{3mm}\nin {\small 2010 Mathematics Subject Classification: 13P10, 13F20, 13C14, 05A18. }

\end{minipage}
\end{center}

\vs{4mm} \section{Introduction }

Throughout the paper,  for a positive integer $n$, let $[n]=\{1,2,\ldots,n\}$. For a set $A$, let $A_d$ be the set  of  the subsets of $A$ with cardinality $d$. In particular, for a simplicial complex $\Delta$, let $\Delta_d$ be the set of faces of $\Delta$ with dimension $d-1$. Let $S=K[x_1, \ldots, x_n]$ be the polynomial ring over a field $k$, and let $I$ be a monomial ideal of $S$. Denote by $sm(S)$ and $sm(I)$ the set of square-free monomials in $S$ and $I$ respectively. There is a natural bijection between $sm(S)$ and $2^{[n]}$, denoted by
$$\sigma: \, x_{i_1}x_{i_2}\cdots x_{i_k} \mapsto \{i_1, i_2, \ldots, i_k\}.$$

For other concepts and notations, see references \cite{AM, Eisenbud, HH, Villareal,Zariski}.

Given a simplicial complex $\Delta$, one can define a Stanley-Reisner ideal $I_{\Delta}$ and a facet ideal $I(\Delta)$ corresponding to $\Delta$. Conversely, given a square-free monomial ideal $I$ of $S=K[x_1, \ldots, x_n]$, there are a pair of simplicial complexes related to $I$. One is the {\it facet} complex of $I$ and is denoted by $\delta_{\mathcal{F}}(I)$. $\delta_{\mathcal{F}}(I)$ is generated by the set $\sigma(G(I))$, i.e., $\sigma(G(I))=\{ \sigma(g) \,|\, g \in G(I)\}$ is the set of facets of $\delta_{\mathcal{F}}(I)$, where $G(I)$ is the minimal generating set  of the monomial ideal $I$ of $S$.  The other one is the Stanley-Reisner complex $\delta_{\mathcal{N}}(I)$ of $I$, or alternatively, the non-face complex  of $I$. Note that the Stanley-Reisner ideal of $\delta_{\mathcal{N}}(I)$ is $I$, in other words, $\delta_{\mathcal{N}}(I) = \{ \sigma(g) \,|\, g \in sm(S) \setminus sm(I)\}$.
The above correspondences construct a bridge between algebraic properties of ideals and combinatorial properties of simplicial complexes.
In order to study algebraic properties such as linear resolution of square-free monomial ideals, one usually takes advantage of the structures of  simplicial complexes corresponding to the ideal, see references \cite{Faridi,HHZ,Zheng, CF}.

Throughout the paper, a monomial ideal $I$ is called {\it of degree $d$}  (or alternatively, {\it homogeneous of degree $d$}), if all monomials in  $G(I)$ have the same degree $d$. Note that the {\it degree of a monomial ideal} $I$, denoted by $deg(I)$,  is the maximal degree of the monomials in $G(I)$. Note  the difference between the two phrases.

Recall that a square-free monomial ideal $I$ is called an {\it $f$-ideal}, if both $\delta_{\mathcal{F}}(I)$ and $\delta_{\mathcal{N}}(I)$ have the same $f$-vector. Note that the $f$-vector of a complex $\delta_{\mathcal{N}}(I)$  is essential in the computation of  the Hilbert series of $S/I$, and  in general the $f$-vector of $\delta_{\mathcal{N}}(I)$ is not easy to calculate. Since the correspondence of the complex $\delta_{\mathcal{F}}(I)$ and the ideal $I$ is direct and clear, it is more easier to calculate the $f$-vector of $\delta_{\mathcal{F}}(I)$. So, it is easy to calculate the Hilbert series and study other corresponding properties of $S/I$ while $I$ is an $f$-ideal.

It seems that the original impetus for combining the simplicial complex $\delta_{\mathcal{F}}(I)$ with $\delta_{\mathcal{N}}(I)$ comes from Remark 2 of \cite{Faridi}, while the formal definition of an $f$-ideal first appeared in \cite{AAAB}, in which the authors studied the properties of $f$-ideals of degree 2, and presented an interesting characterization of such ideals. In  \cite{AMBZ}, the authors generalized the characterization for $f$-ideals of degree $d$ ($d\ge 2$),  though their main result seems to be a little bit inaccurate, see Example \ref{$f$-Ideal not unmixed} in this paper. The importance of $f$-ideals of degree 2 lies in the fact that they are unmixed, see Proposition 5.2 of this paper.
In this paper, we determine all $f$-ideals of degree 2, thus providing a class of unmixed monomial ideals.

In this paper, we focus on the following questions:

$(1)$ How to characterize $f$-ideals of degree $d$ directly?

(2) How many $f$-ideals of degree $d$ are there in the polynomial ring $S=K[x_1, \ldots, x_n]$?

(3) Is there any $f$-ideal which is not unmixed?

(4) What can one say about $f$-ideals in general case?

In section  2, we give an answer to questions (1). We give a complete answer to question (2) in sections  3 and 4,  in the case $d=2$ . In section 5, we present a class of $f$-ideals which are not unmixed, and prove further that in the case $d=2$, an ideal is $f$-ideal if and only if it is an unmixed $f$-ideal. Finally, in section 6, we give a preliminary answer to question $(4)$.

In order to compare with the definition of the degree of a monomial ideal, we need the following:

\begin{mdef}\label{ldeg} For a monomial ideal $I$, the minimal degree of monomials in $G(I)$ is called the {\it lower degree} of $I$, denoted by $ldeg(I)$.
\end{mdef}

\begin{lem}\label{facet leq non-face} (\cite{AMBZ} Lemma 3.6) Let $I$ be a square-free monomial ideal of degree $d$. Then for each $0\le i < d-1$, $\delta_{\mathcal{F}}(I)_{i+1} \seq \delta_{\mathcal{N}}(I)_{i+1}$ holds. In particular, $f_i(\delta_{\mathcal{F}}(I)) \leq f_i(\delta_{\mathcal{N}}(I))$ holds for each $0\le i < d-1$.
\end{lem}

\vs{3mm}The proof to the following lemma is direct to check, so we omit its verification:

\begin{lem}\label{lower than ldeg} Let $I$ be a square-free monomial ideal of $S= K[x_1, \ldots, x_n]$ with $ldeg(I) = k$. Then for each $0 < i < k$, $f_{i-1}(\delta_{\mathcal{N}}(I)) = C_n^i$ holds. Furthermore, if $I$ is an $f$-ideal, then $$f_{i-1}(\delta_{\mathcal{F}}(I)) = f_{i-1}(\delta_{\mathcal{N}}(I)) = C_n^i$$ holds for each $0 < i < k$.
\end{lem}

\vs{3mm}The following corollary follows directly from Lemma \ref{lower than ldeg}:

\begin{cor}\label{facet equal non-face}(\cite{AMBZ} Lemma 3.7) If $I$ is an $f$-ideal of degree $d$, then
$$f_{i-1}(\delta_{\mathcal{F}}(I))= f_{i-1}(\delta_{\mathcal{N}}(I)) = C_n^i$$
holds  for each $0 < i < d$.
\end{cor}

\section{ Perfect sets and $f$-ideals of degree $d$ }

In order to characterize $f$-ideals clearly, we need the following definitions.

\begin{mdef}\label{upper lower} For a set of square-free monomials $A$ in $K[x_1,\ldots,x_n]$, the {\it upper generated set} $\sqcup(A)$ of $A$ is defined by
$$\sqcup(A) = \{gx_i \,|\, g\in A, x_i \nmid g, \,1\le i\le n\}.$$
Dually, the {\it lower cover set}\, $\sqcap(A)$ of $A$ is defined by
$$\sqcap(A) = \{h \mid   1 \neq h,\,h = g/x_i \text{ for some $ g \in A $ and some $x_i$ with $x_i \mid g$ } \}.$$
\end{mdef}

\vs{3mm}Similarly, we define $\sqcup^2(A)= \sqcup(\sqcup(A))$, and $\sqcup^{\infty}(A)= \cup_{i=1}^{\infty} \sqcup^i(A)$, $\sqcap^{\infty}(A)= \cup_{i=1}^{\infty} \sqcap^i(A)$. Actually, both $\sqcup^{\infty}(A)$ and $\sqcap^{\infty}(A)$ are finite union.
Denote by $sm(S)_d$ the set of square-free monomials of degree $d$ in $S$.

\begin{mdef}\label{upper lower perfect set} Let $S=K[x_1, \ldots, x_n]$, and let $A \seq sm(S)_d$.
$A$ is called {\it upper perfect}, if $\sqcup(A) = sm(S)_{d+1}$ holds. Dually, $A$ is called {\it lower perfect}, if $\sqcap(A) = sm(S)_{d-1}$ holds. If $A$ is both upper perfect and lower perfect, then $A$ is called $(n, d)^{th}$ {\it perfect }, or alternatively, a {\it perfect subset} of $sm(S)_d$. For a given pair of numbers $(n,d)$, the least number among cardinalities of $(n, d)^{th}$ perfect sets is called the $(n, d)^{th}$ perfect number, and will be denoted by $N_{(n,d)}$.
\end{mdef}


\begin{ex}\label{perfect example} Let $S=K[x_1, x_2, x_3, x_4]$. Consider the following three subsets of $sm(S)_2$:
$$A=\{x_1x_2, x_1x_3, x_1x_4\},\,  B=\{x_1x_2, x_1x_3, x_2x_3\},\, C=\{x_1x_2, x_3x_4\}.$$
It is direct to check that $A$ is lower perfect, $B$ is upper perfect, $C$ is perfect.  Note that $x_2x_3x_4 \not\in \sqcup(A)$, so $A$ is not upper perfect. Since $x_4 \not\in \sqcap(B)$, $B$ is not lower perfect.
\end{ex}

\vs{3mm}With the aid of the bijection $\si: sm(S)\to 2^{[n]}$, we can define an upper generated subset, an lower cover subset and a (lower, upper) perfect subsets of $2^{[n]}$, respectively. For example,  a subset $A$ of $ [n]_d$ is called a perfect set, if $\sigma^{-1}(A)$ is a perfect subset of $sm(K[x_1,\ldots,x_n]_d)$.

\vs{3mm} The following theorem will show how to judge a square-free monomial ideal of degree $d$ to be an $f$-ideal directly and conveniently.

\begin{thm}\label{homogeneous $f$-Ideal} Let $S=K[x_1, \ldots, x_n]$, and let $I$ be a square-free monomial ideal of $S$ of degree $d$ with the minimal generating set $G(I)$. Then $I$ is an $f$-ideal if and only if $\, G(I)$ is $(n, d)^{th}$ perfect  and $|G(I)| = \frac{1}{2}C_n^d$ holds true.
\end{thm}

\begin{prof} If $I$ is an $f$-ideal of degree $d$, then by definition, $\delta_{\mathcal{F}}(I)$ and $\delta_{\mathcal{N}}(I)$ have the same f-vector. In particular, $dim(\delta_{\mathcal{N}}(I))=dim(\delta_{\mathcal{F}}(I)) = d-1$. By the definition of the non-face complex $\delta_{\mathcal{N}}(I)$, $I$ is the Stanley-Reisner ideal of $\delta_{\mathcal{N}}(I)$. Hence $I$ contains every square-free monomial of degree $d+1$, thus $G(I)$ is upper perfect since $G(I)$ is homogeneous of degree $d$. Furthermore, note that every facet of $\delta_{\mathcal{F}}(I)$ has dimension $d-1$, and $f_{d-2}(\delta_{\mathcal{F}}(I))= C_n^{d-1}$ holds by  Corollary \ref{facet equal non-face}, it follows that $G(I)$ is lower perfect, thus $G(I)$ is a perfect subset of $sm(S)_d$. Finally, $|G(I)| = C_n^d/2$ clearly holds true.

Conversely, if $G(I)$ is $(n, d)^{th}$ perfect and $|G(I)| = C_n^d/2$, we claim that $\delta_{\mathcal{N}}(I)$ is the simplicial complex generated by $D = E \cup [n]_{d-1}$, where $E=[n]_d \setminus \sigma(G(I))$. In fact, let $\Delta$ be the simplicial complex generated by $D$. Then clearly the Stanley-Reisner ideal $I_{\Delta}$ of $\Delta$ contains all the monomials in $G(I)$ and thus $I\seq I_{\Delta}$ holds.  Note further that $G(I)$ is upper perfect, it follows that $I$ is the Stanley-Reisner ideal $I_{\Delta}$ of $\Delta$. On the other hand,  each set in $[n]_{d-1}$ is a face of $\delta_{\mathcal{F}}(I)$ since $G(I)$ is lower perfect. Thus $\delta_{\mathcal{F}}(I)$ and $\delta_{\mathcal{N}}(I)$ have the same f-vector, and hence $I$ is an $f$-ideal.
\end{prof}

\vs{3mm}\nin{\bf Remark 2.5.} By Definition \ref{upper lower perfect set}, Theorem \ref{homogeneous $f$-Ideal} actually provides a rather simple algorithm for listing all the perfect subsets of $2^{[n]}$ of degree $d$:

1. List all elements of $[n]_d,\, [n]_{d+1}, \, [n]_{d-1}$ respectively.

2. For each $A\seq [n]_d$, if $|A|=\frac{1}{2}C_n^d$, then go to the next step.

3. List $\sqcup(A)$ and $\sqcap(A)$, and then check if  $\sqcup(A)=[n]_{d+1}$ and $\sqcap(A)=[n]_{d-1}$ holds true. Note that $A$ is perfect if and only if both equalities hold.

\section{Perfect number $N_{(n,2)}$ and existence of $(n, 2)^{th}$ $f$-ideals}

In the following, an ideal $I$ is called  {\it an $(n, d)^{th}$ $f$-ideal} if $I$ is an $f$-ideal of $K[x_1, \ldots, x_n]$ of degree $d$. We denote by $V(n, d)$ the set of all $(n, d)^{th}$ $f$-ideals. By the characterization of Theorem \ref{homogeneous $f$-Ideal}, $G(I)$ is an $(n, d)^{th}$ perfect set if $I$ is an $(n, d)^{th}$ $f$-ideal. We denote by $U(n, d)$ the set of $(n, d)^{th}$  $f$-ideals whose minimal generating set {\bf contains} a least $(n, d)^{th}$ perfect set, where a  least $(n, d)^{th}$ perfect set  is an  $(n, d)^{th}$ perfect set with cardinality  $N_{(n, d)}$.

In this section, we mainly study $V(n, 2)$. It is easy to see that $N_{(3,2)} = 2$ and $V(3, 2) = \emptyset$. So, assume $n \geq 4$ in the following. Note that $2 \nmid C_n^2$ holds whenever $n=4k+2$ or $n=4k+3$, thus in these cases,  $V(n, 2) = \emptyset$. In the following, we only consider the case when $n=4k$ or $n=4k+1$. Clearly, it is  important to find a perfect set with the least cardinality. For this purpose, we begin with the calculation of the perfect number $N_{(n,2)}$ for general $n\ge 4$. Combining with Proposition \ref{perfect graph}, the following theorem can be proved by Turan's theorem. In order to read conveniently, we show a direct proof in the following.





\begin{thm}\label{perfect number of degree 2} Let $k$ be a positive integer, and let $n\ge 4$.  Then the perfect number $N_{(n,2)}$ is given under the following rules:
\begin{equation}\label{eq: perfect number}
N_{(n,2)} = \begin{cases}
k^2-k,  &  \text{if } n=2k; \\
k^2,  &  \text{if } n=2k+1.
\end{cases}
\end{equation}
\end{thm}

\begin{prof} We will prove the conclusion by the following two steps.

Step 1: we want to estimate the lower bound of the cardinalities of $(n,2)^{th}$ perfect sets.

Let $A$ be an $(n,2)^{th}$  perfect set. Denote $A_{i+} = \{j \,|\, j>i$ and $x_ix_j \in A\}$.

Substep 1: If $x_ix_j \in A$ holds for each pair of $i, j \in [n]$, then the cardinality of $A$ is very large. Without loss of generality, assume that $x_1x_2 \not\in A$. Note that $A$ is upper perfect, hence for every $i \in [n] \setminus [2]$, $x_1x_2x_i \in \sqcup(A)$. So, either $x_1x_i \in A$ or $x_2x_i \in A$. Hence $A_{1+} \cup A_{2+} = [n] \setminus [2]$ and thus $|A_{1+}| + |A_{2+}| \geq n-2$ holds.

Substep 2: This substep is similar to step 1. In fact, if $x_ix_j \in A$ holds for each pair of $i, j \in [n] \setminus [2]$, then the cardinality of $A$ is very large. Without loss of generality, assume that $x_3x_4 \not\in A$. Since for every $i \in [n] \setminus [4]$, $x_3x_4x_i \in \sqcup(A)$. So, either $x_3x_i \in A$ or $x_4x_i \in A$. Hence $A_{3+} \cup A_{4+} = [n] \setminus [4]$, thus $|A_{3+}| + |A_{4+}| \geq n-4$ holds.

Continuing the substeps if necessary. It is easy to see that for each positive integer $l$ such that $2l \leq n$, we have $A_{(2l-1)+} \cup A_{2l+} = [n] \setminus [2l]$ and $|A_{(2l-1)+}| + |A_{2l+}| \geq n-2l$ hold. Therefore, $|A| = \sum_{i=1}^n |A_{i+}| \geq (n-2) + (n-4) + \cdots $ holds. Then we proceed the  calculation  in the following two subcases:

If $n=2k$ for some positive integer $k$, then $|A| \geq k^2-k$. If $n=2k+1$ for some positive integer $k$, then $|A| \geq k^2$.

Step 2: We will show that $N_{(n, 2)}$ can get to the lower bound.

Consider also the two subcases: $n=2k$ and $n=2k+1$.
The following discussion are based on the assumption $n=2k$, and the other case is similar to construct, so we omit the details.

Assume $n=2k$ for some positive integer $k$. Based on the discussion of step 1, in order to show that $N_{(n, 2)}$ can get to the lower bound $k^2-k$, it is suffice to show that for each $l=k, k-1,  \ldots, 1$, we can distribute the elements of $[n] \setminus [2l]$ to $A_{(2l-1)+}$ and $A_{2l+}$ properly.

Substep 0: Set $A_{n+} = \emptyset$, and set $A_{(n-1)+} = \emptyset$.

Substep 1: Set $A_{(n-2)+} = \{n\}$, and set $A_{(n-3)+} = \{n-1\}$.

Analysis: By now, we can make sure that $x_ix_jx_t \in  \sqcup(A)$ for $\{i, j, t\} \seq [n] \setminus [n-4]$.

Substep 2: Set $A_{(n-4)+} = \{n-2\} \cup A_{(n-2)+}$, and set $A_{(n-5)+} = \{n-3\} \cup A_{(n-3)+}$.

Analysis: By now, we can make sure that $x_ix_jx_t \in  \sqcup(A)$ for $\{i, j, t\} \seq [n] \setminus [n-6]$.

Substeps 3, 4, et. al. are similar to substeps 1 and 2. In general, for a positive integer $l \leq n/2$, set $A_{2l+} = \{2l+2\} \cup A_{(2l+2)+}$, and set $A_{(2l-1)+} = \{2l+1\} \cup A_{(2l+1)+}$. It is not hard to see that $x_ix_jx_t \in  \sqcup(A)$ for $\{i, j, t\} \seq [n] \setminus [2l-2]$. In fact, if $\{2l-1, 2l\} \cap \{i, j, t\} = \emptyset$, then by the previous substep, the conclusion is true. If $\{2l-1, 2l\} \cap \{i, j, t\} \neq \emptyset$, without loss of generality, assume $i=2l$. Now consider $j$ and $t$:  if one of them is in $A_{2l+}$, then the conclusion is true. In the other case, $\{j, t\} \seq A_{(2l-1)+}$ holds, thus the conclusion is also true since $x_jx_t \in A$ by the definition of $A_{(2l-1)+}$.

Finally, set $A_{2+} = \{4\} \cup A_{4+}$, and set $A_{1+} = \{3\} \cup A_{3+}$. We also have $x_ix_jx_t \in  \sqcup(A)$ for $\{i, j, t\} \seq [n]$.

By now, we get an upper perfect set $A$. It is easy to see that the set $A$ is also lower perfect, and $|A|=k^2-k$ holds.    \end{prof}

\begin{rem}\label{two complete graphs} The proof to Theorem \ref{perfect number of degree 2}  also answers the afore-mentioned question: How to find an $(n, 2)^{th}$ perfect set with the least cardinality? In fact, what is  needed is to decompose the set $[n]$ into a disjoint union of two subsets $B$ and $C$ {\it uniformly}, namely, $||B|-|C|| \leq 1$. Then set $A= \{x_ix_j \,|\, i, j \in B,$ or $i, j \in C\}$. By the above theorem, $A$ is an $(n, 2)^{th}$ perfect set. Actually, it is easy to check directly that $A$ is a perfect set, and the cardinality of $A$ is equal to the $(n, 2)^{th}$ perfect number $N_{(n, 2)}$, which provides another new understanding of the formula in Theorem \ref{perfect number of degree 2}, i.e.,
\begin{equation}\label{eq: perfect number new understanding}
N_{(n,2)} = \begin{cases}
C_k^2+C_k^2=k^2-k,  &  \text{if } n=2k; \\
C_k^2+C_{k+1}^2=k^2,  &  \text{if } n=2k+1.
\end{cases}
\end{equation}
Note that any set $D$ such that $A \seq D \seq sm(S)_2$ is also an $(n,2)^{th}$  perfect set.  It follows clearly from the definition of an $(n,2)^{th}$   perfect set.
\end{rem}

Now we are ready to settle the existence of $(n, 2)^{th}$ $f$-ideals:

\begin{prop}\label{V(n,2) nonempty}  $V(n,2) \neq \emptyset$ if and only if $n=4k$ or $n=4k+1$ for some positive integer $k$.
\end{prop}

\begin{prof} The necessary part is clear. For the sufficient part, it is suffice to show that the $(n, 2)^{th}$ perfect number is not greater than $C_n^2/2$ in the two cases respectively. If $n=4k$, then by Theorem \ref{perfect number of degree 2} $N_{(n, 2)} = 4k^2-2k$ and $C_n^2/2 = 4k^2-k$, so $N_{(n, 2)} <  C_n^2/2$. If $n=4k+1$, then $N_{(n, 2)} = 4k^2$ and $C_n^2/2 = 4k^2+k$, so we also have $N_{(n, 2)} <  C_n^2/2$.
\end{prof}

\section{Structure of $V(n, 2)$}

For  a nontrivial subset $B$ of $[n]$, let $\ol{B}$ be the complement of $B$ in $[n]$ and let
$$W_B = \{x_ix_j \,|\, i, j \in B\,\, or\,\,i, j \in \ol{B}\}$$
 be a subset of $sm(S)$. Clearly $W_B = W_{\ol{B}}$ holds, and $W_B$ is an $(n, 2)^{th}$ perfect set. Such a way to construct a perfect set is called a {\it Two Part Complete Construction}. A subset $A$ of $sm(S)_2$ is called satisfying {\it Two Part Complete Structure}, abbreviated as TPCS, if there exists a $B \seq [n]$, such that $W_B \seq A$. If further $|B|=l$, then $A$ is called satisfying $l^{th}$ TPCS. An $f$-ideal $I$ is called {\it of $l$ type}, if $G(I)$ satisfies $l^{th}$ TPCS. Denote by $W_l$ the set of  $f$-ideals of  $l$ type in $S$. It is easy to see that $U(n, 2) = W_{2k}$ holds, if $n=4k$ or $n=4k+1$ for some positive integer $k$. We begin with a counting formula for $|U(n,2)|$:


\begin{prop}\label{U(n,2) cardinality} If $k$ is a positive integer, then
\begin{equation}\label{eq: U(n,2) number}
|U(n,2)| = \begin{cases}
\frac{1}{2}C_{4k}^{2k}C_{4k^2}^k,  &  \text{if } n=4k; \\
C_{4k+1}^{2k}C_{4k^2+2k}^k,  &  \text{if } n=4k+1; \\
0,  &  \text{otherwise } 
\end{cases}
\end{equation}
\end{prop}

\begin{prof}We only prove the case when $n=4k$, and the other cases are similar to this one. Assume $I \in U(n,2)$, where $n=4k$. Since $U(n,2)=W_{2k}$, there exists a subset $B \seq [n]$ with $|B|=2k$, such that $W_B \seq G(I)$ holds. We claim that such a subset $B$ is unique, i. e., if there exists another $B_1 \seq [n]$ with $|B_1|=2k$  such that $W_{B_1} \seq G(I)$, then $\{B, \ol{B}\}=\{B_1, \ol{B_1}\}$ holds. In fact, note that both $|G(I)|=C_{4k}^2/2=4k^2-k$ and $|W_B|=2C_{2k}^2=4k^2-2k$ hold, hence there are at most $k$ monomials  in $G(I)\setminus W_B$. Now assume to the contrary that $\{B, \ol{B}\} \neq \{B_1, \ol{B_1}\}$ holds, and assume without loss of generality further that $1, 2 \in B$, $1 \in B_1$ and $2 \not\in B_1$ hold. Then $W_{B_1}$ contains half of the monomials in $W= \{x_1x_j \,|\, j \not\in B\} \cup \{x_2x_j \,|\, j \not\in B\}$. Let $M \seq W$ be such that $M \seq W_{B_1}$ and $|M|=2k$.  Then $M \seq G(I)\setminus W_B$ and $|G(I)\setminus W_B| \leq k$ hold, a contradiction. The contradiction shows the uniqueness of the set $\{B,\ol{B}\}$.

In order to count the cardinality of $U(n,2)$, we need first choose a $2k$ set $B$ randomly, then choose $k$ monomials of $sm(S)_2 \setminus W_B$ randomly. Note that $W_B = W_{\ol{B}}$ holds, thus $|U(n,2)| = \frac{1}{2}C_{4k}^{2k}C_{C_{4k}^2-2C_{2k}^2}^k = \frac{1}{2}C_{4k}^{2k}C_{4k^2}^k$ also holds. This completes the proof.
\end{prof}

\vs{3mm}In the rest part of this section, we will consider a possible decomposition of $V(n, 2)$ into a disjoint union of the afore-mentioned $W_l$.
For this purpose, it is natural to ask the following interesting question: Is there any $f$-ideal who is of no $l$ type?

The following example gives an immediate answer to the question. In the following, we will show that this is the only kind of the example.

\begin{ex}\label{5 cycle}Let $S=K[x_1, x_2, x_3, x_4, x_5]$. It is direct to check that
$$I = \langle x_1x_2, x_2x_3, x_3x_4, x_4x_5, x_1x_5 \rangle$$
is an $f$-ideal, but $I$ is not of $l$ type for any $l$.
\end{ex}

\vs{3mm}How to find further $f$-ideals which are not of $l$ type for any $l$? In order to answer this question,  we need a new idea to construct an $f$-ideal.

Let $S=K[x_1, \ldots, x_n]$, and let $\tau$ be a bijection sending a subset $A$ of $sm(S)_2$ to a graph $T$ whose vertices are $v_1, \ldots, v_n$, such that $v_iv_j \in E(T)$ holds if and only if  $x_ix_j \in A$, where $E(T)$ is the edge set of $T$.

The above example shows that if $T$ is a cycle with 5 vertices, then the ideal generated by $\tau^{-1}(T)$ is an $f$-ideal. Such a class of $f$-ideals will be denoted by $C_5$, which consists of 12 $f$-ideals.

\begin{prop}\label{perfect graph}Let $A \seq sm(S)_2$. Then

$(1)$ $A$ is upper perfect if and only if $\omega(\ol{\tau(A)}) \leq 2$ holds, where $\ol{\tau(A)}$ is the complement graph of $\tau(A)$.

$(2)$ $A$ is lower perfect if and only if for each $i \in [n]$, $d(v_i) < n-1$ holds in the graph $\ol{\tau(A)}$.
\end{prop}

\begin{prof} (1) We will prove it by reduction to an absurdity. For the sufficiency part, assume to the contrary that $A$ is not upper perfect.
Then there exists  a subset $\{i, j, t\} \seq [n]$ such that none of $x_ix_j, x_ix_t, x_jx_t$ is  in $A$, hence
$\{v_iv_j, v_iv_t, v_jv_t\}\cap E(\tau(A)) = \emptyset$ holds, thus $\{v_i, v_j, v_t\}$ is a clique in $\ol{\tau(A)}$, contradicting $\omega(\ol{\tau(A)}) \leq 2$.
The necessity part is similar to get, if we reverse the above discussion.

(2) It is not hard to see that, $A$ is lower perfect if and only if there exists no vertex $v_i \in V(\tau(A))$ with $d(v_i)=0$ in $\tau(A)$, and the latter holds if and only if
 for each $i \in [n]$, $d(v_i) < n-1$ holds in the graph $\ol{\tau(A)}$. This completes the proof.
\end{prof}

\vs{3mm}By the above proposition, Theorem 3.1 is equivalent to the classical Turan's theorem, and the proof to Theorem 3.1 is an alternative proof to Turan's theorem.

\begin{prop}\label{(l, n-l) type bipartite graph}If $I$ is an $(n, 2)^{th}$ $f$-ideal, then $I$ is of $l$ type for some $1 \leq l \leq \lfloor n/2 \rfloor$ if and only if $\ol{\tau(G(I))}$ is a bipartite graph.
\end{prop}

\begin{prof} For the necessity part, assume first that $I$ is of $l$ type. By definition, there exists a subset $B \seq [n]$ with $|B|=l$, such that $W_B \seq G(I)$ holds. Hence $K_B\cup K_{\ol{B}} \seq E(\tau(G(I)))$, where $K_B=\{v_iv_j \,|\, i, j \in B\}$ and $K_{\ol{B}}=\{v_iv_j \,|\, i, j \not\in B\}$. It is easy to see that $\ol{\tau(G(I))}$ is a bipartite graph, with two parts corresponding to $B$ and $\ol{B}$, respectively.

Reversing the above discussion, we get the proof of the sufficiency part.
\end{prof}

\vs{3mm}By Proposition \ref{perfect graph} and Proposition \ref{(l, n-l) type bipartite graph}, the following lemma is clear.

\begin{lem}\label{FC}$I$ is an $(n, 2)^{th}$ $f$-ideal which is not of \,$l$ type for any $l$, if and only if $\ol{\tau(G(I))}$ satisfies the following four conditions (abbreviated as FC in what follows):

$(1)$  For each $i \in [n]$, $d(v_i) < n-1$ holds in $\ol{\tau(G(I))}$.

$(2)$\, $\omega(\ol{\tau(G(I))}) =2$.

$(3)$\, $|E(\ol{\tau(G(I))})|=\frac{C_n^2}{2}$.

$(4)$\, $\ol{\tau(G(I))}$ is not a bipartite graph.
\end{lem}

\vs{3mm}Note that a square-free monomial ideal $I$ of degree $2$ is an $f$-ideal if and only if $\ol{\tau(G(I))}$ satisfies the above conditions $(1),(2)$ and $(3)$.

Our idea to deal with the above FC question is illustrated the following:

Based on the condition (4), construct our model by the rules (1) and (2), and then check the model by the condition (3).

Note that a graph is bipartite if and only if the graph contains no odd cycle, so we will construct our model on an odd cycle. We find an amazing result as the following theorem shows.

\begin{thm}\label{only C_5} If $n \neq 5$, then $V(n, 2) = \bigcup_{l=1}^{\lfloor n/2 \rfloor} W_l$, which is a mutually disjoint union of the $W_l$'s.
\end{thm}

\begin{prof} Note that $V(n, 2) = \cup_{l=1}^{\lfloor n/2 \rfloor} W_l$ holds true,  if and only if each $f$-ideal is of $l$ type for some $l$; and the latter holds if and only if,  there is no graph satisfying the FC. We will show that a graph will not satisfy condition $(3)$ if it satisfies conditions $(2)$ and $(4)$, except for the case $n=5$.

Assume that $T$ is a graph satisfying conditions $(2)$ and $(4)$. Since  $T$ is not a bipartite graph, there exists at least an odd cycle in $T$. Assume that $D$ is a minimal odd cycle of $T$, with $|V(D)| = 2i+1$. Note that $\omega(T) =2$, so $i \geq 2$. Denote by $|E(D)|$ the edge number of the subgraph induced on $D$, and denote by $|E(B, C)|$ the number of edges, each of which has end vertices in $B$ and $C$ respectively. It is clear that $$|E(T)| = |E(D)| + |E(T \setminus D)| + |E(D, T \setminus D)|.$$
holds. Note that $|E(D)| = 2i+1$ holds, since $D$ is a minimal cycle. Since there exists no triangles in $T$, it is not hard to see that $$|E(D, T \setminus D)| \leq (n-2i-1)i$$ holds, since $D$ is an odd cycle.  We will discuss $|E(T \setminus D)|$ in the following two subcases:

If $n=2k$ for some positive $k$, then $|V(T \setminus D)| = 2k-2i-1$ holds. It follows from Turan's theorem that $|E(T \setminus D)| \leq (k-i)(k-i-1)$ hold, hence we get
$$|E(T)| = |E(D)| + |E(T \setminus D)| + |E(D, T \setminus D)|$$
 $$\leq (2i+1)+(2k-2i-1)i+(k-i)(k-i-1) = k^2-k-i^2+2i+1.$$
 Note that $C_n^2/2 = k^2 -k/2$, thus $$C_n^2/2 - |E(T)| \geq k/2+i^2-2i-1 = k/2+(i-1)^2-2$$ holds. Since $i \geq 2$ and $2k > 2i+1$, $C_n^2/2 - |E(T)| > 0$ holds. This shows that there is no graph satisfying FC when $n=2k$.

If $n=2k+1$,
then $|V(T \setminus D)| = 2k-2i$ holds. Again by Turan's theorem, $|E(T \setminus D)| \leq (k-i)^2$ holds, hence we have
$$|E(T)| = |E(D)| + |E(T \setminus D)| + |E(D, T \setminus D)|$$
$$\leq (2i+1)+(2k-2i)i+(k-i)^2 = k^2-i^2+2i+1.$$
Note that $C_n^2/2 = k^2 +k/2$, thus  $$C_n^2/2 - |E(T)| \geq k/2+i^2-2i-1 = k/2+(i-1)^2-2$$ holds true. Then we have $C_n^2/2 - |E(T)| \geq 0$, since $i \geq 2$ and $k \geq i$ hold by assumption. Note further that the equality holds if and only if $k=i=2$.  Thus in this case, there is no graph satisfying FC except $n=5$. This completes the proof.
\end{prof}

\vs{3mm}By the proof of the above theorem, if an $f$-ideal of degree $2$ is not  of $l$ type for any $l$, then it  must be contained in the set  $C_5$, see Example 4.2.

In order to explain the above theorem more precisely, we need the following proposition.

\begin{prop}\label{W not empty} $(1)$ If $n=4k$ for some positive integer $k$, then $W_{2k-i} \neq \emptyset$ holds if and only if $i^2 \leq k$.

$(2)$ If $n=4k+1$, then $W_{2k-i} \neq \emptyset$ holds if and only if $i^2+i \leq k$.
\end{prop}

\begin{prof}(1) Note that $W_{2k-i} \neq \emptyset$ if and only if $C_{2k-i}^2 + C_{2k+i}^2 \leq C_{4k}^2/2$ holds. By direct calculation, the latter holds if and only if $i^2 \leq k$.

(2) It is similar to (1) to check.
\end{prof}

\vs{3mm}The following refines Theorem 4.6:

\begin{thm}\label{exact characterization of V(n, 2)}Let $k$ be a positive integer. Then the following equalities hold true:
\begin{equation}\label{eq: exact V(n,2)}
V(n,2) = \begin{cases}
\underset{0 \leq i \leq \sqrt{k}}{\cup} W_{2k-i},  &  \text{if } n=4k; \\
\underset{0 \leq i \leq \frac{\sqrt{1+4k}-1}{2}}{\cup} W_{2k-i},  &  \text{if } n=4k+1 (k \neq 1); \\
W_2 \cup C_5,  &  \text{if } n=5; \\
\emptyset,  &  \text{if } n=4k+2 \text{ or } n=4k+3.
\end{cases}
\end{equation}
\end{thm}


\begin{rem}\label{(n, 2)^{th} f-ideal construction} Theorem \ref{exact characterization of V(n, 2)} shows the construction of any $(n, 2)^{th}$ $f$-ideal clearly. In the following, we show the construction of an $(n, 2)^{th}$ $f$-ideal while $n=4k$, the other cases are similar to construct.

$(1)$ Choose a nonempty subset $B \seq [n]$, such that $|B| = i \leq \sqrt{k}$;

$(2)$ Let $t=k-i^2$, and choose a subset $E_t \seq sm(S)_2 \setminus W_B$ such that $|E_t|= t$;

$(3)$ Let $I$ be the ideal with the minimal generating set $G(I) = W_B \cup E_t$.
\end{rem}

The proof of the following proposition is similar to Proposition \ref{U(n,2) cardinality}, so we omit it.

\begin{prop}\label{intersection empty} Let $i, j \in [\lfloor n/2 \rfloor]$. Then the following hold:

$(1)$ If $i \neq j$, then $W_i \cap W_j = \emptyset$;

$(2)$ If $I \in W_i$, then there exists a unique subset $B \seq [n]$ with $|B|=i$, such that $W_B \seq G(I)$.
\end{prop}

Note that in the above proposition, the uniqueness of $B$ refers to the uniqueness of $W_B$.


\vs{3mm} 

By Theorem \ref{exact characterization of V(n, 2)} and Proposition \ref{intersection empty}, the following proposition is direct to check, so we omit the proof.


\begin{prop}\label{cardinality of V(n, 2)}Let $k$ be a positive integer. Then the following formula holds:
\begin{equation}\label{eq: cardinality V(n,2)}
|V(n,2)| = \begin{cases}
\frac{1}{2}C_{4k}^{2k}C_{4k^2}^k + \underset{1 \leq i \leq \sqrt{k}}{\sum} C_{4k}^{2k-i}C_{4k^2-i^2}^{k-i^2},  &  \text{if } n=4k; \\
\underset{0 \leq i \leq \frac{\sqrt{1+4k}-1}{2}}{\sum} C_{4k+1}^{2k-i}C_{4k^2+2k-i-i^2}^{k-i-i^2},  &  \text{if } n=4k+1 (k \neq 1); \\
72,  &  \text{if } n=5; \\
0,  &  \text{if } n=4k+2 \text{ or } n=4k+3.
\end{cases}
\end{equation}
\end{prop}

\vs{3mm} By now, the structure of $V(n, 2)$ is completely characterized. However, a complete characterization of $V(n, d)$ for $d > 2$ is still open.

\section{Unmixed $f$-ideals}

It is known that Cohen-Macaulay property is very important in commutative algebra. In \cite{Faridi}, Faridi proved that a Cohen-Macaulay simplicial complex is unmixed. So, it is essential to study the $f$-ideals which are unmixed. Recall that an ideal $I$ is called unmixed if $codim(I)=codim(P)$ holds for all prime ideals minimal over $I$. Recall also the following famous Unmixed Theorem: If $I$ is generated by $r$ elements and $codim(I)=r$, then $I$ is unmixed (see, e.g., \cite[Corollary 18.14]{Eisenbud}).

The following example shows that an $f$-ideal need not to be unmixed.

\begin{ex}\label{$f$-Ideal not unmixed} Let $S=K[x_1, x_2, x_3, x_4, x_5]$, and let
$$I= \langle x_1x_2x_3, x_1x_2x_4, x_1x_2x_5, x_3x_4x_5, x_2x_3x_4 \rangle.$$
 It is not hard to check that $G(I)$ is perfect and $|G(I)|=5=C_5^2/2$, which satisfies the condition of Theorem \ref{homogeneous $f$-Ideal}. Hence $I$ is an $f$-ideal. But the standard primary decomposition of $I$ is $I= \langle x_2, x_5 \rangle \cap \langle x_2, x_3 \rangle \cap \langle x_2, x_4 \rangle \cap \langle x_1, x_4 \rangle \cap \langle x_1, x_3 \rangle \cap \langle x_3, x_4, x_5 \rangle$, which shows that $I$ is not unmixed.
\end{ex}

\vs{3mm} However, when using formulae of section 3 to consider the $f$-ideals of degree $2$, we rediscover the following surprising property, which constitutes the main part of  \cite[Theorem 3.5]{AAAB}. Note that our approach is combinatoric, and is quite different from the proof of  \cite[Theorem 3.5]{AAAB}.

\begin{prop}\label{$f$-Ideal = unmixed $f$-Ideal}  Let $S=K[x_1, \ldots, x_n]$ and let $I$ be a square-free monomial ideal of $S$ of degree 2. If $I$ is an $f$-ideal, then $I$ is unmixed.
\end{prop}

\begin{prof} 
Assume on the contrary that $I$ is not unmixed. By Corollary 1.11 of \cite{Faridi}, $\delta_{\mathcal{N}}(I)$ is not pure. Assume without loss of generality that $\{n\}$ is a facet of $\delta_{\mathcal{N}}(I)$. Then it is easy to see that each of $\{1, n\}, \{2, n\}, \ldots, \{n-1, n\}$ is a facet of $\delta_{\mathcal{F}}(I)$, hence $G_1 = \{x_1x_n, x_2x_n, \ldots, x_{n-1}x_n\} \seq G(I)$.
Further more, $G(I)$ contains at least another $(n-1, 2)^{th}$ upper perfect set, denoted by $G_2$, since $x_{i_1}x_{i_2}x_{i_3} \in \sqcup(G(I))$ for each $\{i_1, i_2, i_3\} \seq [n-1]$.  Hence the cardinality of $G(I)$ is not less than the sum of the cardinalities of the above two parts $G_1$ and $G_2$.

Now that $I$ is an $f$-ideal, by Proposition \ref{V(n,2) nonempty}, $n=4k$ or $n=4k+1$ holds for some positive integer $k$. In the following, we will make use of the formula in Theorem \ref{perfect number of degree 2} to estimate the cardinality of $G(I)$ in the two cases respectively.

If $n=4k$, then $$ |G(I)| \geq |G_1|+|G_2| \geq (n-1) + N_{(n-1, 2)} = (4k-1) + (2k-1)^2= 4k^2,$$
and $C_n^2/2 = 4k^2-k$, so $|G(I)| > C_n^2/2$. In the case, there is a contradiction to Theorem 2.4.

If $n=4k+1$, then $$ |G(I)| \geq |G_1|+|G_2| \geq (n-1) + N_{(n-1, 2)} = 4k + (4k^2-2k)= 4k^2+2k,$$
and $C_n^2/2 = 4k^2+k$, so $|G(I)| > C_n^2/2$,  another contradiction.

This completes the proof.
\end{prof}

\vs{3mm} It is known that a square-free monomial ideal $I$ is unmixed, if and only if $\delta_{\mathcal{N}}(I)$ is a pure simplicial complex (\cite{Faridi}). So, we have the following proposition:

\begin{prop}\label{unmixed homogeneous $f$-Ideal} Let $S=K[x_1, \ldots, x_n]$. If $I$ is an $f$-ideal of $S$ of degree $d$, then $I$ is unmixed if and only if $sm(S)_d \setminus G(I)$ is lower perfect.
\end{prop}

\section{$F$-ideals in general case}

For a square-free monomial ideal $I$, denote $G(I) = \cup_{i=1}^k G_{d_i}$, in which $G_{d_i}$
consists of the generators of degree $d_i$.  As Theorem \ref{homogeneous $f$-Ideal}, the following proposition is direct to check, and we omit the verification.

\begin{prop}\label{nonhomogeneous $f$-Ideal} Let $I$ be a square-free monomial ideal of $S=K[x_1, \ldots, x_n]$, with the minimal generating set $G(I) = \cup_{i=1}^k G_{d_i}$. Then $I$ is an $f$-ideal if and only if the following conditions hold:

$(1)$ For each positive $l \in \{d_1, \ldots, d_k\}$, $$|G_l|= \frac{1}{2}(C_n^l- |\cup_{d_i >  l}(\sqcap^{d_i-l}(G_{d_i}))|- |\cup_{d_i <  l}(\sqcup^{l-d_i}(G_{d_i}))|\,).$$

$(2)$ For each positive $l \not\in \{d_1, \ldots, d_k\}$, $$\underset{d_i >  l}{\cup} \sqcap^{d_i-l}(G_{d_i})  = sm(S)_l \setminus \underset{d_i < l}{\cup}\sqcup^{l-d_i}(G_{d_i}).$$
\end{prop}

\vs{3mm} Actually, we have another way to understand the above proposition directly and clearly. We need the following lemma, which can be checked directly and thus we omit the proof.

\begin{lem}\label{inclusion lemma}If $I$ is a square-free monomial ideal, then

$(1)$ $\sigma(\sqcap^{\infty}(G(I))) = \delta_{\mathcal{F}}(I) \cap \delta_{\mathcal{N}}(I)$;

$(2)$ $\sigma(\sqcup^{\infty}(G(I))) \cap \delta_{\mathcal{F}}(I) = \emptyset$ and $\sigma(\sqcup^{\infty}(G(I))) \cap \delta_{\mathcal{N}}(I) = \emptyset$.
\end{lem}

\begin{thm}\label{general $f$-Ideal} Let $I$ be a square-free monomial ideal of $S=K[x_1, \ldots, x_n]$, with the minimal generating set $G(I) = \cup_{i=1}^k G_{d_i}$. Then $I$ is an $f$-ideal if and only if $$|G_l|= \frac{1}{2}(C_n^l- |\cup_{d_i >  l}(\sqcap^{d_i-l}(G_{d_i}))|- |\cup_{d_i <  l}(\sqcup^{l-d_i}(G_{d_i}))|\,)$$  holds for each $l \in [n]$.
\end{thm}

\begin{prof} For each $l \in [n]$, denote by $A_l$ the faces in $(\delta_{\mathcal{N}}(I) \setminus \delta_{\mathcal{F}}(I)) \cap [n]_l$. Note that $sm(S)_l$ is a disjoint union of four parts: $$sm(S)_l = G_l \cup \cup_{d_i >  l}(\sqcap^{d_i-l}(G_{d_i})) \cup \cup_{d_i <  l}(\sqcup^{l-d_i}(G_{d_i})) \cup \sigma^{-1}(A_l).$$ By Lemma \ref{inclusion lemma}, $f_{l-1}(\delta_{\mathcal{F}}(I)) = |G_l|+|\cup_{d_i >  l}(\sqcap^{d_i-l}(G_{d_i}))|$ and $f_{l-1}(\delta_{\mathcal{N}}(I)) = |\cup_{d_i >  l}(\sqcap^{d_i-l}(G_{d_i}))|+|\sigma^{-1}(A_l)|$. Thus $I$ is an $f$-ideal if and only if $f_{l-1}(\delta_{\mathcal{F}}(I)) = f_{l-1}(\delta_{\mathcal{N}}(I))$ holds for each $l$, and the latter holds if and only if $|G_l| = |\sigma^{-1}(A_l)| =\frac{1}{2}(C_n^l- |\cup_{d_i >  l}(\sqcap^{d_i-l}(G_{d_i}))|- |\cup_{d_i <  l}(\sqcup^{l-d_i}(G_{d_i}))|\,)$  holds for each $l$.
\end{prof}

\vs{3mm} Note that in Proposition \ref{nonhomogeneous $f$-Ideal}, if $l \not\in \{d_1, \ldots, d_k\}$, then $|G_l|=0$. It is not hard to see that Proposition \ref{nonhomogeneous $f$-Ideal} and Theorem \ref{general $f$-Ideal} are equivalent.

Even though the abstract properties of $f$-ideals in general case are characterized in Proposition \ref{nonhomogeneous $f$-Ideal} and Theorem \ref{general $f$-Ideal}, it is still not easy to show an example of an $f$-ideal  which is  not homogeneous of degree $d$ for any $d\ge 2$.

\begin{cor}\label{$f$-Ideal perfect} Let $I$ be an $f$-ideal of $S=K[x_1, \ldots, x_n]$. For an integer $l \in [n]$ with  $l < ldeg(I)$,  $sm(S)_l \seq \sqcap^{\infty}(G(I))$ holds true. On the other hand, if $l > deg(I)$, then $sm(S)_l \seq \sqcup^{\infty}(G(I))$ holds.
\end{cor}

\vs{3mm} It is easy to see that Theorem \ref{homogeneous $f$-Ideal} is a special case of Theorem \ref{general $f$-Ideal}.

\end{document}